\documentclass[numreferences]{kluwer}
\usepackage{amssymb,graphicx}

\newcommand{\ds}{\displaystyle}
\newtheorem{defn}{Definition}

\newcommand{\ason}{\renewcommand{\arraystretch}{1.5}}
\newcommand{\asoff}{\renewcommand{\arraystretch}{1.0}}

\begin{document}
\begin{article}
\begin{opening}

\title
{A Family of Runge-Kutta Methods with Zero Phase-Lag and Derivatives for the Numerical Solution of
the Schr\"{o}dinger Equation and Related Problems}

\author{Z.A. \surname{Anastassi}\thanks{e-mail: zackanas@uop.gr}}
\author{D.S. \surname{Vlachos}\thanks{e-mail: dvlachos@uop.gr}}
\author{T.E. \surname{Simos}\thanks{Highly Cited Researcher, Active Member of the European Academy of Sciences and Arts, Address: Dr. T.E. Simos, 26 Menelaou Street, Amfithea - Paleon Faliron, GR-175 64 Athens, GREECE, Tel: 0030 210 94 20 091, e-mail: tsimos.conf@gmail.com, tsimos@mail.ariadne-t.gr}}
 \institute{Laboratory of Computer Sciences, Department of Computer Science and Technology, Faculty of Sciences and Technology, University of
Peloponnese, GR-22 100 Tripolis, GREECE}

\runningtitle{Runge-Kutta Methods with Zero Phase-Lag and First Derivative}
\runningauthor{Z.A. Anastassi, D.S Vlachos, T.E. Simos}

\begin{abstract}
We construct a family of two new optimized explicit Runge-Kutta methods with zero phase-lag and derivatives for the numerical solution of the time-independent radial Schr\"odinger equation and related ordinary differential equations with oscillating solutions. The numerical results show the superiority of the new technique of nullifying both the phase-lag and its derivatives.
\end{abstract}


\classification{PACS}{0.260, 95.10.E}
\end{opening}

\section{Introduction}
\label{Intro}

Much research has been done on the numerical integration of the radial Schr\"{o}dinger equation:

\begin{equation}
\label{Schrodinger}
    y''(x) = \left( \frac{l(l+1)}{x^{2}}+V(x)-E \right) y(x)
\end{equation}

\noindent where $\frac{l(l+1)}{x^{2}}$ is the \textit{centrifugal potential}, $V(x)$ is the \textit{potential}, $E$ is
the \textit{energy} and $W(x) = \frac{l(l+1)}{x^{2}} + V(x)$ is the \textit{effective potential}. It is valid that
${\mathop {\lim} \limits_{x \to \infty}} V(x) = 0$ and therefore ${\mathop {\lim} \limits_{x \to \infty}} W(x) = 0$.

Many problems in chemistry, physics, physical chemistry, chemical physics, electronics etc., are expressed by equation
(\ref{Schrodinger}).

In this paper we will study the case of $E>0$. We divide $[0,\infty]$ into subintervals $[a_{i},b_{i}]$ so that $W(x)$
is a constant with value ${\mathop{W_{i}}\limits^{\_}}$. After this the problem (\ref{Schrodinger}) can be expressed by
the approximation

\begin{equation}
\begin{array}{l}
\label{Schrodinger_simplified}
y''_{i} = ({\mathop{W}\limits^{\_}} - E)\,y_{i}, \quad\quad \mbox{whose solution is}\\
y_{i}(x) = A_{i}\,\exp{\left(\sqrt{{\mathop{W}\limits^{\_}}-E}\,x\right)} +
B_{i}\,\exp{\left(-\sqrt{{\mathop{W}\limits^{\_}}-E}\,x\right)}, \\
 A_{i},\,B_{i}\,\in {\mathbb R}.
\end{array}
\end{equation}

There has been an extended bibliography on the development and analysis of numerical methods for the efficient solution of the Schr\"odinger equation: see for example \cite{royal}-\cite{jnaiam3_11}.

\section{Basic theory}
\label{theory}
\subsection{Explicit Runge-Kutta methods}
\label{theory_rk}

An $s$-stage explicit Runge-Kutta method used for the computation of the approximation of $y_{n+1}(x)$, when $y_{n}(x)$
is known, can be expressed by the following relations:

\begin{eqnarray}
\label{RK_gen}
\nonumber   y_{n + 1} = y_{n}+{\sum\limits_{i = 1}^{s} {b_{i}}}\,k_{i}\\
            k_{i} = h\, f\left(x_{n} + c_{i} h,\,y_{n} + h\,{\sum\limits_{j = 1}^{i - 1} {a_{ij}\, k_{j}} } \right),\;
            i = 1,\ldots,s
\end{eqnarray}

\noindent where in this case $f\left(x,y(x)\right) = \left( W(x) - E \right) \, y(x)$.

Actually to solve the second order ODE (\ref{Schrodinger}) using first order numerical method (\ref{RK_gen}),
(\ref{Schrodinger}) becomes:

\begin{equation}
\begin{array}{l}
\label{Schrodinger2}
   z'(x) = \left( W(x) - E \right) \, y(x)\\
   y'(x) = z(x)
\end{array}
\end{equation}

\noindent while we use two sets of equations (\ref{RK_gen}): one for $y_{n + 1}$ and one for $z_{n + 1}$.

\noindent The method shown above can also be presented using the Butcher table below:

\begin{equation}
\begin{array}{c|ccccc}
\label{table_RK}
0\\
c_{2} & a_{21}\\
c_{3} & a_{31} & a_{32}\\
\vdots& \vdots& \vdots\\
c_{s} & a_{s1}& a_{s2}& \ldots & a_{s,s-1}& \\
\hline
      & b_{1} & b_{2} & \ldots & b_{s-1}  & b_{s}\\
\end{array}
\end{equation}

\noindent Coefficients $c_{2}$, \ldots, $c_{s}$ must satisfy the equations:

\vspace{-10pt}
\begin{equation}
\label{eq_tree2}
    c_{i} = {\sum\limits_{j = 1}^{i-1} {a_{ij},\;i = 2,\ldots,s}}
\end{equation}

\begin{defn}
\label{defn_tree5} \emph{ \cite{butcher}
    A Runge-Kutta method has algebraic order $p$ when the method's series
    expansion agrees with the Taylor series expansion in the $p$ first terms:}
    $y^{( {n} )}( {x} ) = y_{app.}^{( {n} )} ( {x} )$, $\;\; n=1,2,\ldots,p.$
\end{defn}
\par

A convenient way to obtain a certain algebraic order is to satisfy a number of equations derived from Tree Theory.
These equations will be shown during the construction of the new methods.

\subsection{Phase-Lag Analysis of Runge-Kutta Methods}

The phase-lag analysis of Runge-Kutta methods is based on the test equation

\begin{equation}\label{eq_phase_test1}
y' = I \omega y, \quad \omega \in R
\end{equation}

\noindent Application of the Runge-Kutta method described in (\ref{RK_gen}) to the scalar test
equation (\ref{eq_phase_test1}) produces the numerical solution:

\begin{equation}\label{eq_phase_test2}
{y_{n+1}=a^{n}_{*}}y_{n},\;\; a_{*}=A_{s}(v^{2})+ivB_{s}(v^{2}),
\end{equation}

\noindent where $v=\omega h$ and $A_{s}, B_{s}$ are polynomials in $v^{2}$ completely defined by
Runge-Kutta parameters $a_{i,j}$, $b_{i}$ and $c_{i}$, as shown in (\ref{table_RK}).

\begin{defn}\label{defn_dissip}\emph{\cite{royal}}
\emph{ In the explicit \emph{s}-stage Runge-Kutta method, presented in (\ref{table_RK}), the
quantities}
\begin{center}
$t(v)=v-\arg[a_{*}(v)]$, \quad $a(v)=1-|a_{*}(v)|$\\
\end{center}
\emph{are respectively called the \emph{phase-lag} or \emph{dispersion error} and the
\emph{dissipative error}. If $t(v)=O(v^{q+1})$ and $a(v)=O(v^{r+1})$ then the method is said to be
of dispersive order \emph{q} and dissipative order \emph{r}}.
\end{defn}

\section{Construction of the new trigonometrically fitted Runge-Kutta methods}
\label{Construction}

We consider the explicit Runge-Kutta method with 3 stages and 3rd algebraic order given in table (\ref{table_classical}).

\ason
\begin{equation}
\begin{array}{c|cccccccc}
\label{table_classical}
\frac{1}{2}  & \frac{1}{2}\\
1            & -   1            & 2\\
\hline
        & \frac{1}{6}  & \frac{2}{3} & \frac{1}{6}
\end{array}
\end{equation}
\asoff

\noindent We will construct two new optimized methods.

\subsection{First optimized method with zero phase-lag}
\label{Constr1}

In order to develop the new optimized method, we set free $b_3$, while all other coefficients are borrowed from the classical method. We want the phase-lag of the method to be null, so we satisfy the equation $PL=0$, while solving for $b_3$, where

$$PL = 1/6\, \left( 6+ \left( -2-6\,b_{{3}} \right) {v}^{2} \right) \tan  \left( v \right) +{v}^{3}b_{{3}}+1/6\, \left( -5-6\,b_{{3}} \right) v$$

So $b_3$ becomes
$$\ds b_3 = -{\frac {-6\,\tan \left( v \right) +2\,\tan \left( v \right) {v}^{2}+5\,v}{6 v \left( v\tan \left( v \right) -{v}^{2}+1 \right) }}$$
and its Taylor series expansion is
$$\ds b_3 = \frac{1}{6}-\frac{1}{30}\,{v}^{4}-{\frac {4}{315}}\,{v}^{6}+{\frac {17}{2835}}\,{v}^{8}+{\frac {206}{31185}}\,{v}^{10}+{\frac {7951}{12162150}}\,{v}^{12}-\ldots$$

\noindent where  $v=\omega h$, $\omega$ is a real number and indicates the dominant frequency of the problem and $h$ is the step-length of integration.

\subsection{Second optimized method with zero phase-lag and derivative}
\label{Constr1}

As for the development of the second optimized method, we set free $b_2$ and $b_3$, while all other coefficients are borrowed from the classical method. We want the phase-lag and its first derivative of the method to be null, so we satisfy the equations $\{PL=0, PL'=0\}$, while solving for $b_2$ and $b_3$, where

\begin{equation}
\begin{array}{l}
PL = 1/6\, \left( 6+ \left( -3\,b_{{2}}-6\,b_{{3}} \right) {v}^{2} \right) \tan \left( v \right) +v \left( -1/6-b_{{2}}-b_{{3}}+b_{{3}}{v}^{2}  \right)\\
PL' = -v\tan \left( v \right) b_{{2}}-2\,v\tan \left( v \right) b_{{3}}+5/6+\left( \tan \left( v \right)  \right) ^{2}-1/2\,{v}^{2}b_{{2}}\\
-1/2\,{v}^{2}b_{{2}} \left( \tan \left( v \right)  \right) ^{2}+2\,b_{{3}}{v}^{2}-b_{{3}}{v}^{2} \left( \tan \left( v \right)  \right) ^{2}-b_{{2}}
-b_{{3}}
\end{array}
\end{equation}

Then we have

\ason
\begin{equation}
\begin{array}{l}
b_2=\frac{1}{6}\,{\frac {12
\,v+{v}^{3}+\tan \left( v \right) {v}^{2}-12\,\tan \left( v \right) +{
v}^{3} \left( \tan \left( v \right)  \right) ^{2}}{{v}^{2} \left( -3\,
v+\tan \left( v \right) +v \left( \tan \left( v \right)  \right) ^{2}-
\tan \left( v \right) {v}^{2}+{v}^{3}+{v}^{3} \left( \tan \left( v
 \right)  \right) ^{2} \right) }}\\
b_3=\frac{1}{3}\,{\frac {5\,{v}^{3} \left( \tan \left( v \right)
 \right) ^{2}+7\,{v}^{3}-19\,\tan \left( v \right) {v}^{2}+6\,v
 \left( \tan \left( v \right)  \right) ^{2}-6\,v+6\,\tan \left( v
 \right) }{{v}^{2} \left( -3\,v+\tan \left( v \right) +v \left( \tan
 \left( v \right)  \right) ^{2}-\tan \left( v \right) {v}^{2}+{v}^{3}+
{v}^{3} \left( \tan \left( v \right)  \right) ^{2} \right) }}
\end{array}
\end{equation}
\asoff

The Taylor series expansion of the coefficients are given below:

\ason
\begin{equation}
\begin{array}{l}
b_2 = \frac{2}{3}-\frac{2}{15}\,{v}^{2}-{\frac {52}{315}}\,{v}^{4}-{\frac {3526}{14175}}\,{v
}^{6}-{\frac {173788}{467775}}\,{v}^{8}-{\frac {354768808}{638512875}}
\,{v}^{10}-\ldots\\
b_3 = \frac{1}{6}+\frac{2}{15}\,{v}^{2}+{\frac {25}{126}}\,{v}^{4}+{\frac {4201}{14175}}\,{v
}^{6}+{\frac {207349}{467775}}\,{v}^{8}+{\frac {423287713}{638512875}}
\,{v}^{10}+\ldots
\end{array}
\end{equation}
\asoff

\noindent where  $v=\omega h$, $\omega$ is a real number and indicates the dominant frequency of the problem and $h$ is the step-length of integration.

\section{Algebraic order of the new methods}
\label{Order}

The following 4 equations must be satisfied so that the new methods maintain the third algebraic order of the
corresponding classical method (\ref{table_classical}). The number of stages is symbolized by $s$, where $s=4$. Then we are
presenting the Taylor series expansions of the remainders of these equations, that is the difference of the right part
minus the left part.

\begin{equation}
\begin{array}{cc}\label{alg6_gen}
 \textbf{1st Alg. Order (1 equation)}                                           &    \textbf{3rd Alg. Order (4 equations)}\\
 {\sum\limits_{i = 1}^{s} {b_{i}} }  = 1                                        &    {\sum\limits_{i = 1}^{s} {b_{i}} } c_{i} ^{2} = \frac{1}{3}\\
 \textbf{2nd Alg. Order (2 equations)}                                          &   {\sum\limits_{i,j = 1}^{s} {b_{i}} } a_{ij} c_{j} = \frac{1}{6}\\
 {\sum\limits_{i = 1}^{s} {b_{i}} } c_{i} = \frac{1}{2}
\end{array}
\end{equation}

\subsection{Equations remainders for the first method}

We are presenting ${\it Rem}$ which is the remainder for all four equations for the first method:

\begin{equation}
\begin{array}{l}
{\it Rem} = -\frac{1}{30}\,{v}^{4}-{\frac {4}{315}}\,{v}^{6}+{\frac {17}{2835}}\,{v}^{8}+{\frac {206}{31185}}\,{v}^{10}+\ldots
\end{array}
\end{equation}

\subsection{Equations remainders for the second method}

The four remainders of the equations for the second method are:

\ason
\begin{equation}
\begin{array}{l}
{\it Rem}_{{1}} = \frac{1}{30}\,{v}^{4}+\frac{1}{21}\,{v}^{6}+{\frac {113}{1575}}\,{v}^{8}+{\frac {7171}{66825}}\,{v}^{10}+\ldots\\
{\it Rem}_{{2}} = \frac{1}{15}\,{v}^{2}+{\frac {73}{630}}\,{v}^{4}+{\frac {2438}{14175}}\,{v}^{6}+{\frac {24091}{93555}}\,{v}^{8}+{\frac {245903309}{638512875}}\,{v}^{10}+\ldots\\
{\it Rem}_{{3}} = \frac{1}{10}\,{v}^{2}+{\frac {11}{70}}\,{v}^{4}+{\frac {2213}{9450}}\,{v}^{6}+{\frac {54634}{155925}}\,{v}^{8}+{\frac {37177279}{70945875}}\,{v}^{10}+\ldots\\
{\it Rem}_{{4}} = \frac{2}{15}\,{v}^{2}+{\frac {25}{126}}\,{v}^{4}+{\frac {4201}{14175}}\,{v}^{6}+{\frac {207349}{467775}}\,{v}^{8}+{\frac {423287713}{638512875}}\,{v}^{10}+\ldots\\
\end{array}
\end{equation}
\asoff

We see that the two optimized methods retain the third algebraic order, since the constant term of all
the remainders is zero.

\section{Numerical results}
\label{Numerical_results}

\subsection{The inverse resonance problem}
The efficiency of the two new constructed methods will be measured through the integration of problem
(\ref{Schrodinger}) with $l=0$ at the interval $[0,15]$ using the well known Woods-Saxon potential

\begin{eqnarray}
\label{Woods_Saxon} V(x) = \frac{u_{0}}{1+q} + \frac{u_{1}\,q}{(1+q)^2}, \quad\quad q =
\exp{\left(\frac{x-x_{0}}{a}\right)}, \quad
\mbox{where}\\
\nonumber u_{0}=-50, \quad a=0.6, \quad x_{0}=7 \quad \mbox{and} \quad u_{1}=-\frac{u_{0}}{a}
\end{eqnarray}

\noindent and with boundary condition $y(0)=0$.

\noindent The potential $V(x)$ decays more quickly than $\frac{l\,(l+1)}{x^2}$, so for large $x$ (asymptotic region)
the Schr\"{o}dinger equation (\ref{Schrodinger}) becomes

\begin{equation}
\label{Schrodinger_reduced}
    y''(x) = \left( \frac{l(l+1)}{x^{2}}-E \right) y(x)
\end{equation}

\noindent The last equation has two linearly independent solutions $k\,x\,j_{l}(k\,x)$ and $k\,x\,n_{l}(k\,x)$, where
$j_{l}$ and $n_{l}$ are the \textit{spherical Bessel} and \textit{Neumann} functions. When $x \rightarrow \infty$ the
solution takes the asymptotic form

\begin{equation}
\label{asymptotic_solution}
\begin{array}{l}
 y(x) \approx A\,k\,x\,j_{l}(k\,x) - B\,k\,x\,n_{l}(k\,x) \\
\approx D[sin(k\,x - \pi\,l/2) + \tan(\delta_{l})\,\cos{(k\,x - \pi\,l/2)}],
\end{array}
\end{equation}

\noindent where $\delta_{l}$ is called \textit{scattering phase shift} and it is given by the following expression:

\begin{equation}
\tan{(\delta_{l})} = \frac{y(x_{i})\,S(x_{i+1}) - y(x_{i+1})\,S(x_{i})} {y(x_{i+1})\,C(x_{i}) - y(x_{i})\,C(x_{i+1})},
\end{equation}

\noindent where $S(x)=k\,x\,j_{l}(k\,x)$, $C(x)=k\,x\,n_{l}(k\,x)$ and $x_{i}<x_{i+1}$ and both belong to the asymptotic
region. Given the energy we approximate the phase shift, the accurate value of which is $\pi/2$ for the above problem.

We will use three different values for the energy: i) $989.701916$, ii) $341.495874$ and iii) $163.215341$. As for the
frequency $\omega$ we will use the suggestion of Ixaru and Rizea \cite{ix_ri}:

\begin{equation}
\omega = \cases{ \sqrt{E-50} & $x\in[0,\,6.5]$ \cr
            \sqrt{E}    &$x\in[6.5,\,15]$ \cr}
\end{equation}

\subsection{Nonlinear Problem}
\hspace{12pt} $y''=-100\, y+\sin(y),\;$ with $\; y(0)=0,\: y'(0)=1,\: t\in[0,20\,\pi]$, $y(20\pi)=$ $3.92823991 \,\cdot\, 10 ^{-4}$
and $\omega=10$ as frequency of this problem.

\subsection{Comparison}

\begin{figure}[tbp]
    \includegraphics[width=\textwidth]{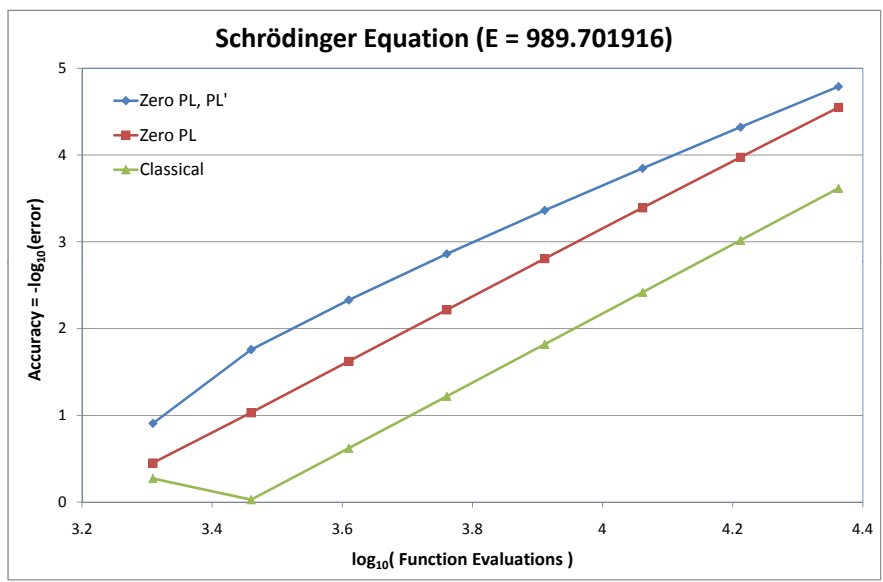}
    \caption{Efficiency for the Schr\"odinger equation using E = 989.701916}
    \label{fig_resonance_989}
\end{figure}

\begin{figure}[tbp]
    \includegraphics[width=\textwidth]{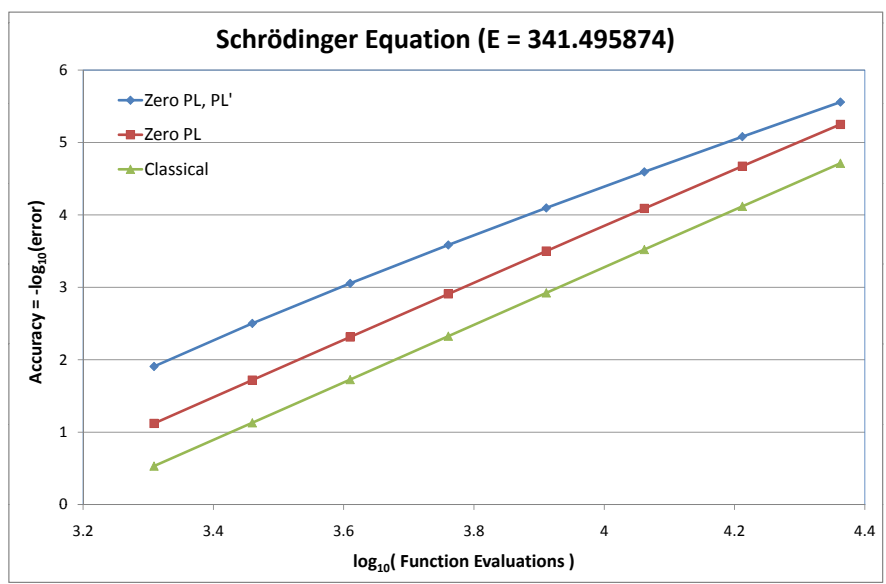}
    \caption{Efficiency for the Schr\"odinger equation using E = 341.495874}
    \label{fig_resonance_341}
\end{figure}

\begin{figure}[tbp]
    \includegraphics[width=\textwidth]{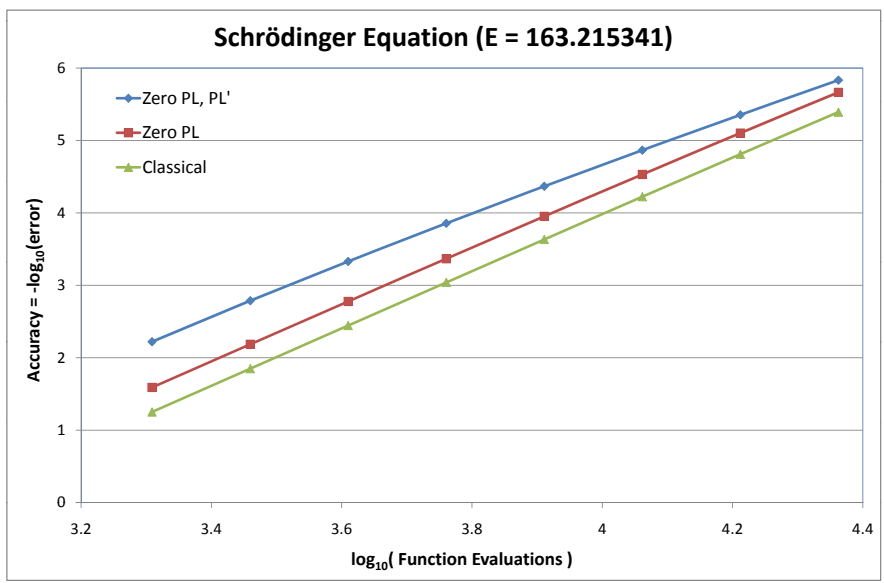}
    \caption{Efficiency for the Schr\"odinger equation using E = 163.215341}
    \label{fig_resonance_163}
\end{figure}

We present the \textbf{accuracy} of the tested methods expressed by the $-\log_{10}$(error at the end point) when
comparing the phase shift to the actual value $\pi/2$ versus the $\log_{10}$(total function evaluations). The
\textbf{function evaluations} per step are equal to the number of stages of the method multiplied by two that is the
dimension of the vector of the functions integrated for the Schr\"odinger ($y(x)$ and $z(x)$). In Figure
\ref{fig_resonance_989} we use $E = 989.701916$, in Figure \ref{fig_resonance_341} $E = 341.495874$ and in Figure
\ref{fig_resonance_163} $E =
163.215341$.

\begin{figure}[tbp]
    \includegraphics[width=\textwidth]{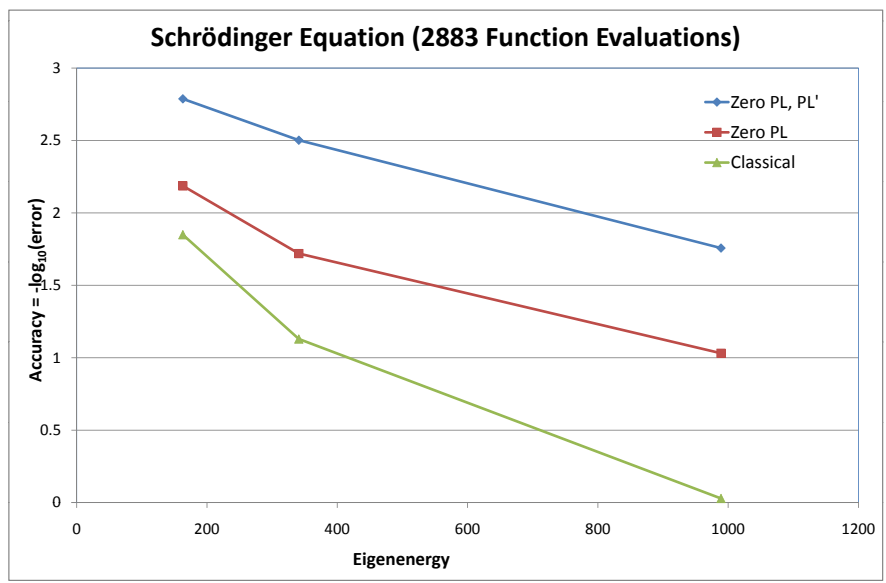}
    \caption{Efficiency for different eigenenergies}
    \label{fig_energies}
\end{figure}

\begin{figure}[tbp]
    \includegraphics[width=\textwidth]{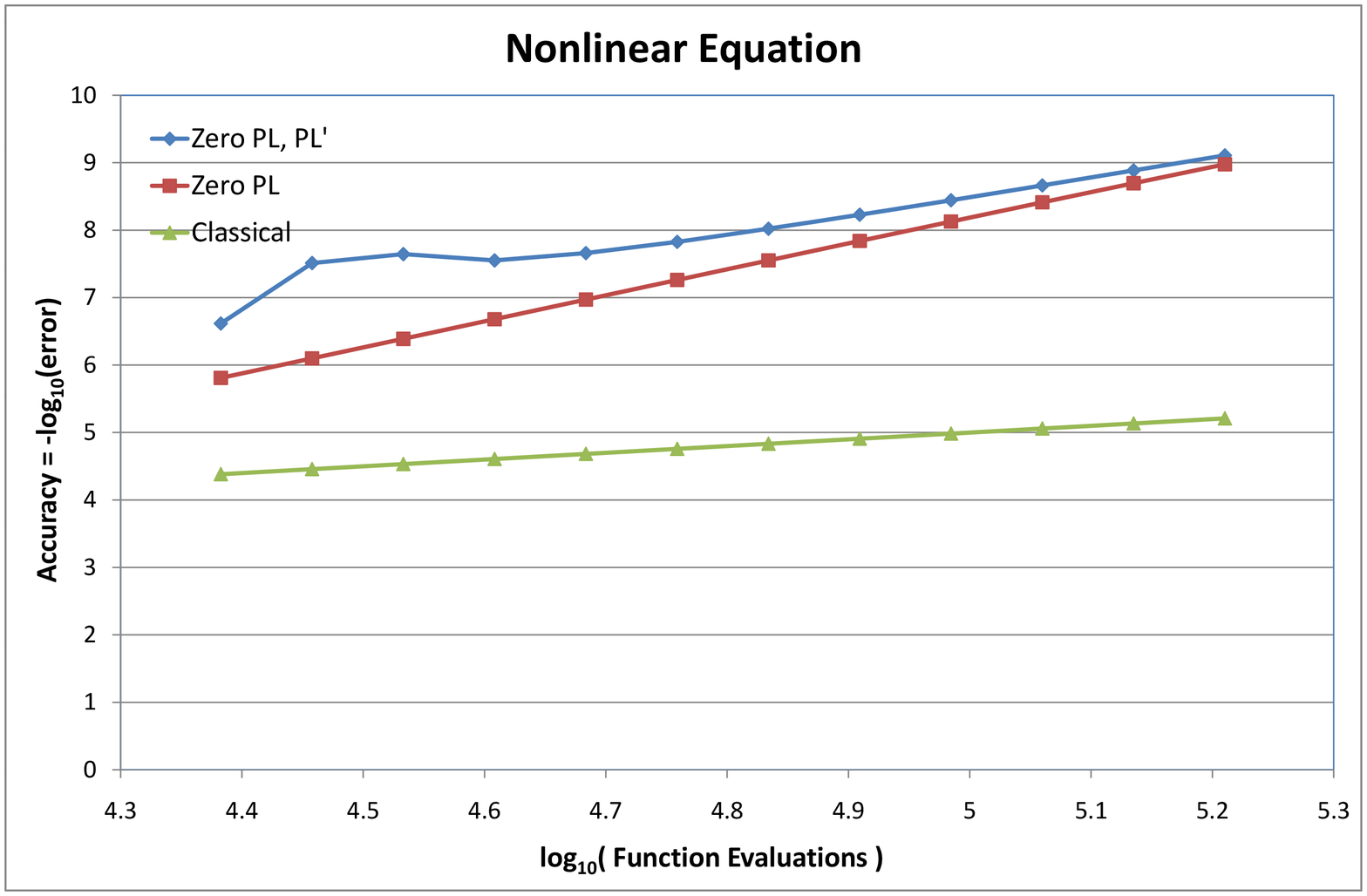}
    \caption{Efficiency for the Nonlinear Problem}
    \label{fig_nonlinear}
\end{figure}

\section{Conclusions}

We compare the two optimized methods and the corresponding classical explicit Runge-Kutta method for the integration of the Schr\"odinger equation and the Nonlinear problem. We see that the second method with the phase-lag and its first derivative nullified is the most efficient in all cases, followed in terms of efficiency by the optimized method with zero phase-lag and then by the corresponding classical method.

\end{article}

\begin{thebibliography}{120}

\bibitem{royal}{Simos T.E., Chemical Modelling - Applications and Theory Vol.1, Specialist Periodical Reports, 2000, The Royal Society of Chemistry, Cambridge, 32-140}

\bibitem{ix_ri} Ixaru L. Gr., Rizea M., A Numerov-like scheme for the numerical solution of the Schr\"{o}dinger equation in the deep continuum spectrum of energies, Comp. Phys. Comm. 19, 23-27 (1980)

\bibitem{butcher} J.C. Butcher, Numerical methods for ordingary differential equations, Wiley (2003)

\bibitem{ix78} L.Gr. Ixaru and M. Micu,
{\it Topics in Theoretical Physics}. Central Institute of Physics,
Bucharest, 1978.

\bibitem{landau} L.D. Landau and F.M.
Lifshitz: {\it Quantum  Mechanics}.  Pergamon,  New York, 1965.

\bibitem{prigogine} I. Prigogine,
Stuart Rice (Eds): Advances in Chemical Physics Vol. 93: New
Methods in Computational Quantum Mechanics, John Wiley \& Sons,
1997.

\bibitem{hertz} G. Herzberg, {\it Spectra of
Diatomic Molecules}, Van Nostrand, Toronto, 1950.

\bibitem{simos00_r} T.E. Simos, Atomic Structure
Computations in Chemical Modelling: Applications and Theory
(Editor: A. Hinchliffe, UMIST), {\it The Royal Society of
Chemistry} 38-142(2000).

\bibitem{simos02_r}  T.E. Simos, Numerical
methods for 1D, 2D and 3D differential equations arising in
chemical problems, {\em Chemical Modelling: Application and
Theory}, The Royal Society of Chemistry, 2(2002),170-270.

\bibitem{simos_wilr} T.E. Simos and
P.S. Williams, On finite difference methods for the solution of
the Schr\"odinger equation, {\it Computers \& Chemistry} {\bf 23}
513-554(1999).

\bibitem{simos90} T.E. Simos: {\it Numerical
Solution of  Ordinary Differential Equations with Periodical
Solution}. Doctoral Dissertation, National Technical University of
Athens, Greece, 1990 (in Greek).

\bibitem{kosim01} A. Konguetsof and T.E.
Simos, On the Construction of exponentially-fitted methods for the
numerical solution of the Schr\"odinger Equation, {\it Journal of
Computational Methods in Sciences and Engineering} {\bf 1}
143-165(2001).

\bibitem{ra78} A.D. Raptis and A.C.
Allison: Exponential - fitting  methods  for the numerical
solution of the Schr\"odinger   equation, \textit{Computer Physics
Communications}, {\bf 14} 1-5(1978).

\bibitem{ra84gr} A.D. Raptis, Exponential
multistep methods for ordinary differential equations, {\it Bull.
Greek Math. Soc.} {\bf 25} 113-126(1984).

\bibitem{ix84} L.Gr.  Ixaru, Numerical Methods
for Differential Equations and Applications, Reidel, Dordrecht -
Boston - Lancaster, 1984.

\bibitem{wilsim02} T. E. Simos, P. S.
Williams: A New Runge-Kutta-Nystrom Method with Phase-Lag of Order
Infinity for the Numerical Solution of the Schr\"odinger Equation,
{\it MATCH Commun. Math. Comput. Chem.} {\bf 45} 123-137(2002).

\bibitem{simos03} T. E. Simos, Multiderivative
Methods for the Numerical Solution of the Schr\"odinger Equation,
{\it MATCH Commun. Math. Comput. Chem.} {\bf 45} 7-26(2004).

\bibitem{ra83} A.D. Raptis, Exponentially-fitted solutions of the
eigenvalue Shr\"odinger equation with automatic error control,
{\it Computer Physics Communications}, {\bf 28} 427-431(1983)

\bibitem{ra81} A.D. Raptis, On the numerical solution of the Schrodinger
equation, {\it Computer Physics Communications}, {\bf 24}
1-4(1981)

\bibitem{siex99} Zacharoula Kalogiratou and T.E. Simos, A P-stable exponentially-fitted method for the numerical
integration of the Schr\"odinger equation, {\it Applied
Mathematics and Computation}, {\bf 112} 99-112(2000).

\bibitem{lambert} J.D. Lambert and I.A.
Watson, Symmetric multistep methods for periodic initial values
problems, {\it J. Inst. Math. Appl.} {\bf 18} 189-202(1976).

\bibitem{rapsim91} A.D. Raptis and T.E. Simos, A four-step phase-fitted method for the numerical  integration
of  second  order  initial-value  problem, {\it BIT}, {\bf 31}
160-168(1991).

\bibitem{henrici} Peter Henrici, {\it Discrete variable methods in ordinary differential equations}, John Wiley
\& Sons, 1962.

\bibitem{chawla83} M.M. Chawla, Uncoditionally stable Noumerov-type methods for second order differential equations, {\it BIT},
{\bf 23} 541-542(1983).

\bibitem{chawla84} M. M. Chawla and P. S. Rao, A Noumerov-type method with minimal phase-lag for the integration of second order periodic initial-value problems, {\it Journal of Computational and Applied Mathematics}  {\bf 11(3)}
277-281(1984)

\bibitem{berghe2004} Liviu Gr. Ixaru and Guido Vanden Berghe,
Exponential Fitting, Series on Mathematics and its Applications,
Vol. 568, Kluwer Academic Publisher, The Netherlands, 2004.

\bibitem{ixaru85} L. Gr. Ixaru and M. Rizea, Comparison of some four-step methods for the numerical solution of the Schr\"odinger equation, {\it Computer
Physics Communications}, {\bf 38(3)} 329-337(1985)

\bibitem{simos1}
Z.A. Anastassi, T.E. Simos, A family of exponentially-fitted Runge-Kutta methods with exponential order up to three for the numerical solution of the Schr\"odinger equation, {\it J. Math. Chem} {\bf 41 (1)} 79-100 (2007)

\bibitem{simos2}
T. Monovasilis, Z. Kalogiratou , T.E. Simos, Trigonometrically
fitted and exponentially fitted symplectic methods for the
numerical integration of the

Schr\"odinger equation, {\it J. Math. Chem} {\bf 40 (3)} 257-267
(2006)

\bibitem{simos3}
G. Psihoyios, T.E. Simos, The numerical solution of the radial
Schr\"odinger equation via a trigonometrically fitted family of
seventh algebraic order

Predictor-Corrector methods, {\it J. Math. Chem} {\bf 40 (3)}
269-293 (2006)

\bibitem{simos4}
T.E. Simos, A four-step exponentially fitted method for the
numerical solution of the Schr\"odinger equation, {\it J. Math.
Chem} {\bf 40 (3)} 305-318

(2006)

\bibitem{simos5}
T. Monovasilis, Z. Kalogiratou , T.E. Simos, Exponentially fitted
symplectic methods for the numerical integration of the
Schr\"odinger equation {\it J. Math. Chem} {\bf 37 (3)} 263-270
(2005)

\bibitem{simos6}
Z. Kalogiratou , T. Monovasilis, T.E. Simos, Numerical solution of
the two-dimensional time independent Schr\"odinger equation with
Numerov-type methods {\it J. Math. Chem} {\bf 37 (3)} 271-279
(2005)

\bibitem{simos7}
Z.A. Anastassi, T.E. Simos, Trigonometrically fitted Runge-Kutta
methods for the numerical solution of the Schr\"odinger equation
{\it J. Math. Chem} {\bf 37 (3)} 281-293 (2005)

\bibitem{simos8}
G. Psihoyios, T.E. Simos, Sixth algebraic order trigonometrically
fitted predictor-corrector methods for the numerical solution of
the radial Schr\"odinger

equation, {\it J. Math. Chem} {\bf 37 (3)} 295-316 (2005)

\bibitem{simos9}
D.P. Sakas, T.E. Simos, A family of multiderivative methods for
the numerical solution of the Schr\"odinger equation, {\it J.
Math. Chem} {\bf 37 (3)}

317-331 (2005)

\bibitem{simos10}
T.E. Simos, Exponentially - fitted multiderivative methods for the
numerical solution of the Schr\"odinger equation, {\it J. Math.
Chem} {\bf 36 (1)} 13-27 (2004)

\bibitem{simos11}
K. Tselios, T.E. Simos, Symplectic methods of fifth order for the
numerical solution of the radial Shrodinger equation, {\it J.
Math. Chem} {\bf 35 (1)}

55-63 (2004)

\bibitem{simos12}
T.E. Simos, A family of trigonometrically-fitted symmetric methods
for the efficient solution of the Schr\"odinger equation and
related problems {\it J. Math. Chem} {\bf 34 (1-2)} 39-58 JUL 2003

\bibitem{simos13}
K. Tselios, T.E. Simos, Symplectic methods for the numerical
solution of the radial Shr\"odinger equation, {\it J. Math. Chem}
{\bf 34 (1-2)} 83-94 (2003)

\bibitem{simos14}
J. Vigo-Aguiar J, T.E. Simos, Family of twelve steps exponential
fitting symmetric multistep methods for the numerical solution of
the Schr\"odinger equation, {\it J. Math. Chem} {\bf 32 (3)}
257-270 (2002)

\bibitem{simos15}
G. Avdelas, E. Kefalidis, T.E. Simos, New P-stable eighth
algebraic order exponentially-fitted methods for the numerical
integration of the Schr\"odinger equation, {\it J. Math. Chem}
{\bf 31 (4)} 371-404 (2002)

\bibitem{simos16}
T.E. Simos, J. Vigo-Aguiar, Symmetric eighth algebraic order
methods with minimal phase-lag for the numerical solution of the
Schr\"odinger equation {\it J. Math. Chem} {\bf 31 (2)} 135-144
(2002)

\bibitem{simos17}
Z. Kalogiratou , T.E. Simos, Construction of trigonometrically and
exponentially fitted Runge-Kutta-Nystrom methods for the numerical
solution of the Schr\"odinger equation and related problems a
method of 8th algebraic order, {\it J. Math. Chem} {\bf 31 (2)}
211-232

\bibitem{simos18}
T.E. Simos, J. Vigo-Aguiar, A modified phase-fitted Runge-Kutta
method for the numerical solution of the Schr\"odinger equation,
{\it J. Math. Chem} {\bf 30 (1)} 121-131 (2001)

\bibitem{simos19}
G. Avdelas, A. Konguetsof, T.E. Simos, A generator and an
optimized generator of high-order hybrid explicit methods for the
numerical solution of the Schr\"odinger equation. Part 1.
Development of the basic method, {\it J. Math. Chem} {\bf 29 (4)}
281-291 (2001)

\bibitem{simos20}
G. Avdelas, A. Konguetsof, T.E. Simos, A generator and an
optimized generator of high-order hybrid explicit methods for the
numerical solution of the Schr\"odinger equation. Part 2.
Development of the generator; optimization of the generator and
numerical results, {\it J. Math. Chem} {\bf 29 (4)} 293-305 (2001)

\bibitem{simos21}
J. Vigo-Aguiar, T.E. Simos, A family of P-stable eighth algebraic
order methods with exponential fitting facilities, {\it J. Math.
Chem} {\bf 29 (3)} 177-189 (2001)

\bibitem{simos22}
T.E. Simos, A new explicit Bessel and Neumann fitted eighth
algebraic order method for the numerical solution of the
Schr\"odinger equation {\it J. Math. Chem} {\bf 27 (4)} 343-356
(2000)

\bibitem{simos23}
G. Avdelas, T.E. Simos, Embedded eighth order methods for the
numerical solution of the Schr\"odinger equation, {\it J. Math.
Chem} {\bf 26 (4)} 327-341 1999,

\bibitem{simos24}
T.E. Simos, A family of P-stable exponentially-fitted methods for
the numerical solution of the Schr\"odinger equation, {\it J.
Math. Chem} {\bf 25 (1)} 65-84 (1999)

\bibitem{simos25}
T.E. Simos, Some embedded modified Runge-Kutta methods for the
numerical solution of some specific Schr\"odinger equations, {\it
J. Math. Chem} {\bf 24 (1-3)} 23-37 (1998)

\bibitem{simos26}
T.E. Simos, Eighth order methods with minimal phase-lag for
accurate computations for the elastic scattering phase-shift
problem, {\it J. Math. Chem} {\bf 21 (4)} 359-372 (1997)

\bibitem{jnaiam1}
P. Amodio, I. Gladwell and G. Romanazzi, Numerical Solution of
General Bordered ABD Linear Systems by Cyclic Reduction, {\it
JNAIAM J. Numer. Anal. Indust. Appl. Math} {\bf 1(1)} 5-12(2006)

\bibitem{jnaiam2}
S. D. Capper, J. R. Cash and D. R. Moore, Lobatto-Obrechkoff
Formulae for 2nd Order Two-Point Boundary Value Problems, {\it
JNAIAM J. Numer. Anal. Indust. Appl. Math} {\bf 1(1)} 13-25 (2006)

\bibitem{jnaiam3}
S. D. Capper and D. R. Moore, On High Order MIRK Schemes and
Hermite-Birkhoff Interpolants, {\it JNAIAM J. Numer. Anal. Indust.
Appl. Math} {\bf 1(1)} 27-47 (2006)

\bibitem{jnaiam4}
J. R. Cash, N. Sumarti, T. J. Abdulla and I. Vieira, The
Derivation of Interpolants for Nonlinear Two-Point Boundary Value
Problems, {\it JNAIAM J. Numer. Anal. Indust. Appl. Math} {\bf
1(1)} 49-58 (2006)

\bibitem{jnaiam5}
J. R. Cash and S. Girdlestone, Variable Step Runge-Kutta-Nyström
Methods for the Numerical Solution of Reversible Systems,  {\it
JNAIAM J. Numer. Anal. Indust. Appl. Math} {\bf 1(1)} 59-80 (2006)

\bibitem{jnaiam6}
Jeff R. Cash and Francesca Mazzia, Hybrid Mesh Selection
Algorithms Based on Conditioning for Two-Point Boundary Value
Problems, {\it JNAIAM J. Numer. Anal. Indust. Appl. Math} {\bf
1(1)} 81-90 (2006)

\bibitem{jnaiam7}
Felice Iavernaro, Francesca Mazzia and Donato Trigiante, Stability
and Conditioning in Numerical Analysis, {\it JNAIAM J. Numer.
Anal. Indust. Appl. Math} {\bf 1(1)} 91-112 (2006)

\bibitem{jnaiam8}
Felice Iavernaro and Donato Trigiante, Discrete Conservative
Vector Fields Induced by the Trapezoidal Method, {\it JNAIAM J.
Numer. Anal. Indust. Appl. Math} {\bf 1(1)} 113-130 (2006)

\bibitem{jnaiam9}
Francesca Mazzia, Alessandra Sestini and Donato Trigiante, BS
Linear Multistep Methods on Non-uniform Meshes, {\it JNAIAM J.
Numer. Anal. Indust. Appl. Math} {\bf 1(1)} 131-144 (2006)

\bibitem{jnaiam10}
L.F. Shampine, P.H. Muir, H. Xu, A User-Friendly Fortran BVP
Solver, {\it JNAIAM J. Numer. Anal. Indust. Appl. Math} {\bf 1(2)}
201-217 (2006)

\bibitem{jnaiam11}
G. Vanden Berghe and M. Van Daele, Exponentially- fitted
Störmer/Verlet methods, {\it JNAIAM J. Numer. Anal. Indust. Appl.
Math} {\bf 1(3)} 241-255 (2006)

\bibitem{jnaiam12}
L. Aceto, R. Pandolfi, D. Trigiante, Stability Analysis of Linear
Multistep Methods via Polynomial Type Variation, {\it JNAIAM J.
Numer. Anal. Indust. Appl. Math} {\bf 2(1-2)} 1-9 (2007)

\bibitem{psih} G. Psihoyios, A Block Implicit Advanced Step-point (BIAS) Algorithm for Stiff Differential
Systems, {\it Computing Letters} {\bf 2(1-2)} 51-58(2006)

\bibitem{enright} W.H. Enright, On the use of 'arc length' and 'defect'
for mesh selection for differential equations, {\it Computing
Letters} {\bf 1(2)} 47-52(2005)

\bibitem{simosnew1}
T.E. Simos, P-stable Four-Step Exponentially-Fitted Method for the
Numerical Integration of the Schr\"{o}dinger Equation, {\em
Computing Letter} \textbf{1(1)} 37-45(2005).

\bibitem{simos2007}
T.E. Simos, Stabilization of a Four-Step Exponentially-Fitted
Method and its Application to the Schr\"odinger Equation, {\em
International Journal of Modern Physics C} \textbf{18(3)}
315-328(2007).

\bibitem{wang2005}
Zhongcheng Wang, P-stable linear symmetric multistep methods for
periodic initial-value problems, {\em Computer Physics
Communications} \textbf{171} 162–174(2005)

\bibitem{simoscma93} T.E. Simos, A Runge-Kutta Fehlberg method with phase-lag of order
infinity for initial value problems with oscillating solution, {\it
Computers and Mathematics with Applications} {\bf 25} 95-101(1993).
\bibitem{simoscma93b} T.E. Simos, Runge-Kutta interpolants with minimal phase-lag,
{\it Computers and Mathematics with Applications} {\bf 26}
43-49(1993).
\bibitem{simoscma93c} T.E. Simos, Runge-Kutta-Nystr\"om interpolants for the numerical
integration of special second-order periodic initial-value problems,
{\it Computers and Mathematics with Applications} {\bf 26}
7-15(1993).
\bibitem{simoscma94} T.E. Simos and G.V. Mitsou, A family of four-step exponential fitted
methods for the numerical integration of the radial Schr\"odinger
equation, {\it Computers and Mathematics with Applications} {\bf 28}
41-50(1994).
\bibitem{simoscma95} T.E. Simos and G. Mousadis, A two-step method for the numerical
solution of the radial Schrödinger equation, {\it Computers and
Mathematics with Applications} {\bf 29} 31-37(1995).
\bibitem{simoscma96} G. Avdelas and T.E. Simos, Block Runge-Kutta methods for periodic
initial-value problems, {\it Computers and Mathematics with
Applications} {\bf 31} 69- 83(1996).
\bibitem{simoscma96b} G. Avdelas and T.E. Simos, Embedded methods for the numerical
solution of the Schr\"odinger equation, {\it Computers and
Mathematics with Applications} {\bf 31} 85-102(1996).
\bibitem{simoscma96c}  G. Papakaliatakis and T.E. Simos, A new method for the numerical
solution of fourth order BVP's with oscillating solutions, {\it
Computers and Mathematics with Applications} {\bf 32} 1-6(1996).
\bibitem{simoscma97}  T.E. Simos, An extended Numerov-type method for the numerical
solution of the Schr\"odinger equation, {\it Computers and
Mathematics with Applications} {\bf 33} 67-78(1997).
\bibitem{simoscma98} T.E. Simos, A new hybrid imbedded variable-step procedure for the
numerical integration of the Schr\"odinger equation, {\it Computers
and Mathematics with Applications} {\bf 36} 51-63(1998).
\bibitem{simoscma01}  T.E. Simos, Bessel and Neumann Fitted Methods for the Numerical
Solution of the Schr\"odinger equation, {\it Computers \&
Mathematics with Applications} {\bf 42} 833-847(2001).
\bibitem{simoscma03}  A. Konguetsof and T.E. Simos, An exponentially-fitted and
trigonometrically-fitted method for the numerical solution of
periodic initial-value problems, {\it Computers and Mathematics with
Applications} {\bf 45} 547-554(2003).
\bibitem{simoscam05} Z.A. Anastassi and T.E. Simos, An optimized Runge-Kutta method for the solution of orbital problems, {\it Journal of Computational and Applied Mathematics} {\bf 175(1)} 1-9(2005)
\bibitem{simoscam05a} G. Psihoyios and T.E. Simos, A fourth algebraic order trigonometrically fitted predictor-corrector scheme for IVPs with
oscillating solutions, {\it Journal of Computational and Applied
Mathematics} {\bf 175(1)} 137-147(2005)
\bibitem{simoscam05b} D.P. Sakas and T.E. Simos, Multiderivative methods of eighth algrebraic order with minimal phase-lag for the numerical solution of
the radial Schr\"odinger equation,  Journal of Computational and
Applied Mathematics   {\bf 175(1)} 161-172(2005)
\bibitem{simoscam05c} K. Tselios and T.E. Simos, Runge-Kutta methods with minimal dispersion and dissipation for problems arising from computational acoustics, {\it Journal of Computational and Applied Mathematics} {\bf 175(1)} 173-181(2005)
\bibitem{simoscam03} Z. Kalogiratou and T.E. Simos, Newton-Cotes formulae for long-time integration, {\it Journal of Computational and Applied Mathematics} {\bf 158(1)} 75-82(2003)
\bibitem{simoscam03a} Z. Kalogiratou, T. Monovasilis and T.E. Simos, Symplectic integrators for the numerical solution of the Schr\"odinger equation, {\it Journal of Computational and Applied Mathematics} {\bf 158(1)} 83-92(2003)
\bibitem{simoscam03b}  A. Konguetsof and T.E. Simos, A generator of hybrid symmetric four-step methods for the numerical solution of the Schr\"odinger equation, {\it Journal of Computational and Applied Mathematics} {\bf 158(1)}  93-106(2003)
\bibitem{simoscam03c} G. Psihoyios and T.E. Simos, Trigonometrically fitted predictor-corrector methods for IVPs with oscillating solutions, {\it Journal of Computational and Applied Mathematics} {\bf 158(1)} 135-144(2003)
\bibitem{simoscam02} Ch. Tsitouras and T.E. Simos, Optimized Runge-Kutta pairs for problems with oscillating solutions, {\it Journal of Computational and Applied Mathematics} {\bf 147(2)} 397-409(2002)
\bibitem{simoscam99} T.E. Simos, An exponentially fitted eighth-order method for the numerical solution of the Schr\"odinger equation, {\it Journal of Computational and Applied Mathematics} {\bf 108(1-2)}
177-194(1999)
\bibitem{simoscam98} T.E. Simos, An accurate finite difference method for the numerical solution of the Schr\"odinger equation, {\it Journal of Computational and Applied Mathematics} {\bf 91(1)} 47-61(1998)
\bibitem{simoscam97} R.M. Thomas and T.E. Simos, A family of hybrid exponentially fitted predictor-corrector methods for the
numerical integration of the radial Schr\"odinger equation, {\it
Journal of Computational and Applied Mathematics} {\bf 87(2)}
215-226(1997)

\bibitem{anastassi_p1} Z.A. Anastassi and T.E. Simos: Special Optimized Runge-Kutta methods for IVPs with Oscillating Solutions, International Journal of Modern Physics C, 15, 1-15 (2004)

\bibitem{anastassi_p3} Z.A. Anastassi and T.E. Simos: A Dispersive-Fitted and Dissipative-Fitted Explicit Runge-Kutta method for the Numerical Solution of Orbital Problems, New Astronomy, 10, 31-37 (2004)

\bibitem{anastassi_p6} Z.A. Anastassi and T.E. Simos: A Trigonometrically-Fitted Runge-Kutta Method for the Numerical Solution of Orbital Problems, New Astronomy, 10, 301-309 (2005)

\bibitem{anastassi_p9} T.V. Triantafyllidis, Z.A. Anastassi and T.E. Simos: Two Optimized Runge-Kutta Methods for the Solution of the Schr?dinger Equation, MATCH Commun. Math. Comput. Chem., 60, 3 (2008)

\bibitem{anastassi2} Z.A. Anastassi and T.E. Simos, Trigonometrically Fitted Fifth Order Runge-Kutta Methods for the Numerical Solution of the Schr\"{o}dinger Equation, Mathematical and Computer Modelling, 42 (7-8), 877-886 (2005)

\bibitem{anastassi_simos_p8} Z.A. Anastassi and T.E. Simos: New Trigonometrically Fitted Six-Step Symmetric Methods for the Efficient Solution of the Schr\"odinger Equation, MATCH Commun. Math. Comput. Chem., 60, 3 (2008)

\bibitem{panopoulos_match} G.A. Panopoulos, Z.A. Anastassi and T.E. Simos: Two New Optimized Eight-Step Symmetric Methods for the Efficient Solution of the Schr\"{o}dinger Equation and Related Problems, MATCH Commun. Math. Comput. Chem., 60, 3 (2008)

\bibitem{anastassi_simos_p11} Z.A. Anastassi and T.E. Simos: A Six-Step P-stable Trigonometrically-Fitted Method for the Numerical Integration of the Radial Schr\"odinger Equation, MATCH Commun. Math. Comput. Chem., 60, 3 (2008)

\bibitem{anas_p12} Z.A. Anastassi and T.E. Simos, A family of two-stage two-step methods for the numerical integration of the Schr\"odinger equation and related IVPs with oscillating solution, Journal of Mathematical Chemistry, Article in Press, Corrected Proof

\bibitem{simoscam97a} T.E. Simos and P.S. Williams, A finite-difference method for the numerical solution of the Schr\"odinger equation, {\it Journal of Computational and Applied Mathematics} {\bf 79(2)} 189-205(1997)
\bibitem{simoscam96} G. Avdelas and T.E. Simos, A generator of high-order embedded P-stable methods for the numerical solution of the
Schr\"odinger equation, {\it Journal of Computational and Applied
Mathematics} {\bf 72(2)}  345-358(1996)
\bibitem{simoscam96a} R.M. Thomas, T.E. Simos and G.V. Mitsou, A family of Numerov-type exponentially fitted predictor-corrector methods for the
numerical integration of the radial Schr\"odinger equation, {\it
Journal of Computational and Applied Mathematics} {\bf 67(2)}
255-270(1996)
\bibitem{simoscam95} T.E. Simos, A Family of 4-Step Exponentially Fitted Predictor-Corrector Methods for the Numerical-Integration of
The Schr\"odinger-Equation, {\it Journal of Computational and
Applied Mathematics} {\bf 58(3)}   337-344(1995)
\bibitem{simoscam94} T.E. Simos, An Explicit 4-Step Phase-Fitted Method for the Numerical-Integration of 2nd-order Initial-Value Problems, {\it Journal of Computational and Applied Mathematics} {\bf 55(2)} 125-133(1994)
\bibitem{simoscam94a} T.E. Simos, E. Dimas and A.B. Sideridis, A Runge-Kutta-Nystr\"om Method for the Numerical-Integration of Special 2nd-order
Periodic Initial-Value Problems, {\it Journal of Computational and
Applied Mathematics} {\bf 51(3)} 317-326(1994)
\bibitem{simoscam92} A.B. Sideridis and T.E. Simos, A Low-Order Embedded Runge-Kutta Method for Periodic Initial-Value Problems, {\it Journal of Computational and Applied Mathematics} {\bf 44(2)} 235-244(1992)
\bibitem{simoscam92a} T.E. Simos amd A.D. Raptis, A 4th-order Bessel Fitting Method for the Numerical-Solution of the Schr\"Odinger-Equation, {\it
Journal of Computational and Applied Mathematics} {\bf 43(3)}
313-322(1992)
\bibitem{simoscam92b} T.E. Simos, Explicit 2-Step Methods with Minimal Phase-Lag for the Numerical-Integration of Special 2nd-order Initial-Value
Problems and their Application to the One-Dimensional
Schr\"odinger-Equation, {\it Journal of Computational and Applied
Mathematics} {\bf 39(1)} 89-94(1992)
\bibitem{simoscam90} T.E. Simos, A 4-Step Method for the Numerical-Solution of the Schr\"odinger-Equation, {\it Journal of Computational and Applied Mathematics} {\bf 30(3)} 251-255(1990)
\bibitem{simoscam90a} C.D. Papageorgiou, A.D. Raptis and T.E. Simos, A Method for Computing Phase-Shifts for Scattering, {\it Journal of Computational and Applied
Mathematics} {\bf 29(1)} 61-67(1990)
\bibitem{raptis82} A.D. Raptis, Two-Step Methods for the Numerical Solution  of the Schr\"odinger Equation, {\it Computing} {\bf 28} 373-378(1982).
\bibitem{simosijmpc96} T.E. Simos. A new Numerov-type method for computing eigenvalues and resonances of the radial Schr\"odinger
equation, International Journal of Modern Physics C-Physics and
Computers, {\bf 7(1)} 33-41(1996)
\bibitem{simosijqc95} T.E. Simos, Predictor Corrector Phase-Fitted Methods for Y''=F(X,Y) and an Application to the
Schr\"odinger-Equation, International Journal of Quantum Chemistry,
{\bf 53(5)} 473-483(1995)
\bibitem{simosijcm92} T.E. Simos, Two-step almost P-stable complete
in phase methods for the numerical integration of second order
periodic initial-value problems, {\it Inter. J. Comput. Math.} {\bf
46} 77-85(1992).
\bibitem{jnaiam3_1} R. M. Corless, A. Shakoori, D.A. Aruliah, L. Gonzalez-Vega, Barycentric Hermite Interpolants for Event Location in Initial-Value Problems, {\it JNAIAM J. Numer. Anal. Indust. Appl. Math}, 3, 1-16 (2008)
\bibitem{jnaiam3_2} M. Dewar, Embedding a General-Purpose Numerical Library in an Interactive Environment, {\it JNAIAM J. Numer. Anal. Indust. Appl. Math}, 3, 17-26 (2008)
\bibitem{jnaiam3_3} J. Kierzenka and L.F. Shampine, A BVP Solver that Controls Residual and Error, {\it JNAIAM J. Numer. Anal. Indust. Appl. Math}, 3, 27-41 (2008)
\bibitem{jnaiam3_4} R. Knapp, A Method of Lines Framework in Mathematica, {\it JNAIAM J. Numer. Anal. Indust. Appl. Math}, 3, 43-59 (2008)
\bibitem{jnaiam3_5} N. S. Nedialkov and J. D. Pryce, Solving Differential Algebraic Equations by Taylor Series (III): the DAETS Code, {\it JNAIAM J. Numer. Anal. Indust. Appl. Math}, 3, 61-80 (2008)
\bibitem{jnaiam3_6} R. L. Lipsman, J. E. Osborn, and J. M. Rosenberg, The SCHOL Project at the University of Maryland: Using Mathematical Software in the Teaching of Sophomore Differential Equations, {\it JNAIAM J. Numer. Anal. Indust. Appl. Math}, 3, 81-103 (2008)
\bibitem{jnaiam3_7} M. Sofroniou and G. Spaletta, Extrapolation Methods in Mathematica, {\it JNAIAM J. Numer. Anal. Indust. Appl. Math}, 3, 105-121 (2008)
\bibitem{jnaiam3_8} R. J. Spiteri and Thian-Peng Ter, pythNon: A PSE for the Numerical Solution of Nonlinear Algebraic Equations, {\it JNAIAM J. Numer. Anal. Indust. Appl. Math}, 3, 123-137 (2008)
\bibitem{jnaiam3_9} S.P. Corwin, S. Thompson and S.M. White, Solving ODEs and DDEs with Impulses, {\it JNAIAM J. Numer. Anal. Indust. Appl. Math}, 3, 139-149 (2008)
\bibitem{jnaiam3_10} W. Weckesser, VFGEN: A Code Generation Tool, {\it JNAIAM J. Numer. Anal. Indust. Appl. Math}, 3, 151-165 (2008)
\bibitem{jnaiam3_11} A. Wittkopf, Automatic Code Generation and Optimization in Maple, {\it JNAIAM J. Numer. Anal. Indust. Appl. Math}, 3, 167-180 (2008)

\end{thebibliography}
\end{document}